\documentclass[pdf, 10pt]{article}
\usepackage[square,authoryear]{natbib}
\usepackage{graphicx}
\usepackage{amscd,amsthm,amssymb,amsmath}

\def\E{\mathbb{E}}

\newtheorem{Condition}{Condition}[section]
\newtheorem{Remark}{Remark}[section]
\newtheorem{Example}{Example}[section]
\newtheorem{theorem}{Theorem}[section]

\newtheorem{definition}[theorem]{Definition}

\newcommand{\bfi}{\bfseries\itshape}

\begin{document}

\title{Stochastic Variational Partitioned Runge-Kutta Integrators for Constrained Systems}

\author{Nawaf Bou-Rabee$^*$  \and Houman Owhadi\thanks{Applied \& Computational Mathematics (ACM), Caltech, Pasadena, CA 91125 ({\tt nawaf@acm.caltech.edu}, \tt{owhadi@acm.caltech.edu}).} }

\maketitle

\begin{abstract}
Stochastic variational integrators for constrained, stochastic mechanical systems are developed in this paper.  The main results of the paper are twofold: an equivalence is established between a stochastic Hamilton-Pontryagin (HP) principle in generalized coordinates and constrained coordinates via Lagrange multipliers, and variational partitioned Runge-Kutta (VPRK) integrators are extended to this class of systems.   Among these integrators are first and second-order strongly convergent RATTLE-type integrators.  We prove strong order of accuracy of the methods provided.  The paper also reviews the deterministic treatment of VPRK integrators from the HP viewpoint.  
\end{abstract}

 \section{Introduction}

Since the foundational work of Bismut [1981], the field of stochastic geometric mechanics is emerging in response to the demand for tools to analyze the structure of continuous and discrete mechanical systems with uncertainty \citep{Ha2007,VaCi2006,CiLeVa2008,Bi1981, MiReTr2002, MiReTr2003, LaOr2007a, LaOr2007b,MaWi2007}.   Within this context the goal of this paper is to develop efficient, structure-preserving integrators for long-time simulations of constrained, mechanical systems perturbed by white-noise forces and torques. Our strategy is to employ stochastic variational integrators (SVIs) \citep{BoOw2007a}.

Variational integration theory derives integrators for mechanical systems from discrete variational principles \citep{Ve1988,Ma1992,WeMa1997,MaWe2001}.   The theory includes discrete analogs of the Lagrangian, Noether's theorem, the Euler-Lagrange equations, and the Legendre transform.  Variational integrators can readily incorporate holonomic constraints (e.g., via Lagrange multipliers) and non-conservative effects (via their virtual work)  \citep{WeMa1997, MaWe2001}.    Altogether, this description of mechanics stands as a self-contained theory of mechanics akin to Hamiltonian, Lagrangian or Newtonian mechanics.

Variational integrators are {\em symplectic}, i.e.,  the discrete flow map they define exactly preserves the continuous symplectic 2-form.  Using backward error analysis one can show that symplectic integrators applied to Hamiltonian systems nearly preserve the energy of the continuous mechanical system for exponentially long periods of time and that the modified equations are also Hamiltonian \citep{HaLuWa2006}.  Variational integrators are also distinguished by their ability to compute statistical properties of mechanical systems, such as in computing Poincar\'e sections, the instantaneous temperature of a system, etc.

Stochastic variational integrators are an extension of variational integrators to so-called stochastic mechanical systems.  These systems are simple mechanical systems subject to certain random perturbations, and have their origins in Bismut's foundational work \citep{Bi1981}.  Bismut showed that the flow of these stochastic mechanical systems extremize an action integral whose domain is the space of semimartingales on configuration space.  Bismut's work was further enriched and generalized to manifolds by recent work \citep{LaOr2007a, LaOr2007b}.   Lazaro-Cami and Ortega show that this general class of stochastic Hamiltonian systems on manifolds extremizes a stochastic action defined on the space of manifold-valued semimartingales \citep{LaOr2007a}.    Moreover, it has been shown that for a subclass of these systems, one can prove a converse, namely,  a.s.~a curve satisfies so-called stochastic Hamilton's equations if and only if it extremizes a stochastic action \citep{BoOw2007a}.

With this variational principle, one can design SVIs \citep{BoOw2007a}.  Like their deterministic counterparts, these methods have the advantage that they are symplectic, and in the presence of symmetry, satisfy a discrete Noether's theorem.   Moreover, symplectic methods for stochastic mechanical systems on vector spaces have been shown to capture the correct energy behavior even in the presence of dissipation \citep{MiReTr2002, MiReTr2003}.   In particular, the energy injected or dissipated by the symplectic integrator is not an artifical function of the timestep.  Moreover, the energy behavior is accurately captured by the integrator.

These structure-preserving properties are quite crucial.  Consider, for instance, simulation of a simple mechanical system subject to random forces with amplitude $\epsilon$. Suppose further the unperturbed system preserves energy, momentum, and possesses a first integral.  Consider simulating this system with a higher-order accurate method, a standard integrator with simultaneous projection onto energy, momentum, and first integral level sets, and a stochastic variational integrator.  If $\epsilon \ne 0$, a stochastic (physical) perturbation in the energy, momentum, and first integrals will appear.  However, it is not clear how to modify the projection-based method to correctly capture these stochastic perturbations when $\epsilon \ne 0$.  Moreover, the higher-order accurate method requires a time-step smaller than the amplitude of the perturbation in order to accurately represent its effects.  And even then, if the time-span of integration is long enough systematic, artificial drift in these quantities will appear.

On the other hand, the stochastic perturbation in energy and momentum captured by the stochastic variational integrator is mechanical, i.e., it is only due to the $\epsilon$-random forces.  This is because these schemes define flows of discrete mechanical systems.  Moreover, even when the amplitude of perturbation is not small, due to symplecticity we conjecture that SVIs on manifolds not only possess good long-time energy behavior, but also perform well in computing statistical properties, such as  autocorrelation functions and the empirical distribution.   We will investigate these questions in future work.

As far as we can tell, the extension of these structure-preserving integrators to stochastic mechanical systems with holonomic constraints has not been completed.   The main goal of this paper is to extend SVIs to holonomic constraints, and in particular, introduce {\bfi constrained, stochastic variational partitioned Runge-Kutta} (VPRK) methods for such systems.   Within this family we exhibit a first-order strongly convergent, symplectic integrator for constrained mechanical systems.  We use a technique due to Vanden-Eijnden and Cicotti [2006] to prove order of accuracy of the integrators that appear in this paper \citep{VaCi2006}.

Future work will consider extensions of these schemes to Langevin equations with holonomic constraints.  Continuous Langevin processes have been generalized to submanifolds $\Sigma \subset \mathbb{R}^n$ and shown to be ergodic \citep{CiLeVa2008}.  The generalization to constraint submanifolds was done by appending holonomic constraints to Langevin equations.  One can then check that the infinitesimal generator of the constrained Langevin process commutes with the Gibbsian density restricted to $\Sigma$ to determine that the restricted Gibbsian measure is an invariant measure of the constrained Langevin process.  To prove this measure is unique one uses standard arguments based on nondegeneracy of the diffusion and drift vector fields on the momentums.  On this note it would be interesting to ascertain ergodicity of SVIs for constrained Langevin systems.

\section{Constrained VPRK Integrators}

This section is provided to fix notation and clarify some aspects of deterministic constrained VPRK integrators which will be relevant in the stochastic context.  We also provide a novel proof showing that constrained VPRK integrators can be derived directly from a discrete variational principle without explicitly introducing a discrete Legendre transform.  For more details on unconstrained VPRK integrators we refer the reader to  \citep{Bo2007}.

The setting of this section is a real, $n$-dimensional  vector space $Q$ and a mechanical system with smooth, holonomic constraint function, $g: Q \to \mathbb{R}^k$, $k<n$, that has a regular value at $0 \in \mathbb{R}^k$.    The mechanical system's configuration space is given by the {\em constraint submanifold}: $S = g^{-1}(0) \subset Q$.    We will introduce Lagrange multipliers to prove an equivalence between a Hamilton-Pontryagin (HP) variational principle on $S$ and a constrained HP variational principle on $Q$.   We will then show how to discretize this system to obtain variational RATTLE  integrators for constrained mechanical systems.

The variational and symplectic character of VPRK integrators is discussed in \citep{Su1990, MaWe2001, HaLuWa2006}.  In what follows we will explicitly use the HP perspective, and specifically, extend the results in Hairer et al.~[2006] to constrained systems using the HP perspective.  It is worth mentioning, the Hamiltonian counterparts of the constrained VPRK methods, so-called symplectic partitioned Runge-Kutta methods, are also well understood \citep{Ja1996, Re1997, Ha2003}.

\paragraph{Discretization of HP Action}

We will adopt a HP viewpoint to describe the action integral of this mechanical system with constraints.  The HP description unifies the Hamiltonian and Lagrangian descriptions of a mechanical system  \citep{YoMa2006a, YoMa2006b, Bo2007, BoMa2007}.

\begin{definition}
The {\bfi Pontryagin bundle} of a manifold $M$ is defined as $PM = T M \oplus T^* M$.  Fixing the interval $[a,b] \subset \mathbb{R}$ and $x_1, x_2 \in M$, define the {\bfi HP path space} on $PM$ as
\[
\mathcal{C}(PM,x_1,x_2) = \{ (q,v,p) \in C^{\infty}([a,b], PM)~|~q(a)=x_1,~q(b)=x_2 \}  \text{.}
\]
\label{defn:hppathspace}
\end{definition}

The HP path space is a smooth infinite-dimensional manifold.  One can show that its tangent space at  $c = (q,v,p) \in \mathcal{C}( PQ, x_1, x_2)$ consist of maps $w = (q,v,p, \delta q, \delta v, \delta p) \in C^{\infty}([a,b],T(PQ))$  such that  $\delta q(a) = \delta q(b) = 0$.

\begin{definition}
Fix $q_1, q_2 \in N$.  Define the  {\bfi unconstrained HP action} $\mathfrak{G}:  \mathcal{C}(PQ,q_1,q_2) \to \mathbb{R}$  as
\[
\mathfrak{G} =  \int _a^b \left[ L(q, v) dt + \left\langle p, \frac{dq}{dt} - v \right\rangle dt \right] \text{.}
\]
\label{defn:hpaction}
\end{definition}

To discretize $\mathfrak{G}$ we first discretize the kinematic constraint: $dq/dt = v$.  An s-stage RK method is employed to discretize the kinematic constraint.  Let $[a,b]$ and $N$ be given and define the fixed step size $h=(b-a)/N$ and $t_k = h k$, $k=0,...,N$.   The reason for using an s-stage RK discretization of the kinematic constraint is that the theory on such methods (order conditions, stability, and implementation) is mature.  See, for instance, \citep{HaNoWa1993}.

\begin{definition}  
Consider the first order differential equation 
\begin{equation}
\dot{q} = f(t,q),~~~q(0) = q_0,~~~q(t) \in Q \text{.} \label{eq:ivpvs}
\end{equation}
Let $b_i, a_{ij} \in \mathbb{R}$ ($i, j=1,\cdots,s$) and let $c_i = \sum_{j=1}^s a_{ij}$.  An {\bfi s-stage RK} approximation is given by 
\begin{equation} \label{eq:sstagerk}
\begin{cases}
\begin{array}{ccl}
Q_k^i &=& q_k + h \sum_{j=1}^s a_{ij} f(t_k+c_j h, Q_k^j) \text{,} ~~~ i=1, \cdots, s \text{,}   \\
q_{k+1} &=& q_k + h \sum_{j=1}^s b_j f(t_k+c_j h, Q_k^j) \text{.} 
\end{array}
\end{cases}
\end{equation}
The vectors $q_k$ and $Q_k^i$ are called external and internal stage vectors, respectively. 
\end{definition}

It follows that an s-stage RK method is fully determined by its $a$-matrix and $b$-vector which are typically displayed using the so-called Butcher tableau: 
\begin{center}
\begin{tabular}{c|ccc}
    $c_1$ &  $a_{11}$ & $\cdots$ & $a_{1s}$ \\
$\vdots$ & $\vdots$   &             & $\vdots$ \\
  $c_s$  & $a_{s1}$   &  $\cdots$  & $a_{ss}$ \\
\hline  
          & $b_1$       & $\cdots$   & $b_s$ 
\end{tabular}
\end{center}
Suppose that $v(t)$, $t \in [a,b]$, is given.  Then an s-stage RK approximant applied to $\dot{q} = v(t)$ yields:
\begin{equation} \label{eq:sstagerkkc}
\begin{cases}
\begin{array}{ccl}
Q_k^i &= & q_k + h \sum_{j=1}^s a_{ij} v(t_k+c_j h)  \text{,} ~~~ i=1, \cdots, s \text{,} \\
q_{k+1} &=& q_k + h \sum_{j=1}^s b_j v(t_k+c_j h)  \text{,} ~~~ k=0, \cdots, N \text{.} 
\end{array}
\end{cases}
\end{equation}
In what follows $v(t_k+c_i h)$ for $i=1,...,s$ will be introduced as internal stage unknowns and will be determined as a critical point of a discrete action sum.

\paragraph{VPRK Integrator for Constrained Systems}

The VPRK method will be derived from a discretization of the HP action integral in which the kinematic constraint over the kth-time step is replaced with its discrete approximant: (\ref{eq:sstagerkkc}), and the integral of the Lagrangian over the kth-time step is approximated by the following quadrature:
\[
\int_{t_k}^{t_k + h} L(q,v) dt \approx h b_i L(Q_k^i, V_k^i)  \text{.}
\]  
The constraint $g(q)=0$ is enforced for all internal stage positions $\{ Q_k^i \}$ using Lagrange multipliers as follows.
\begin{definition}
Fix two points $q_1$ and $q_2$ on $Q$ and define the {\bfi discrete VPRK path space}  as:
\begin{align*}
\mathcal{C}_d = \{  (q,p, \{Q^i,V^i,P^i \}_{i=1}^s, \{ \Lambda^i \}_{i=1}^s)_d :& \{ t_k \}_{k=0}^N  \to T^*Q \times (TQ \oplus T^*Q)^s \times (\mathbb{R}^k)^s ~|~ \\
&  q(0) = q_1, q(t_N) = q_2 \} \text{,}
\end{align*}
and the {\bfi discrete constrained VPRK action sum} $\mathfrak{G}_d: \mathcal{C}_d \to \mathbb{R}$ by:
\begin{align*}
\mathfrak{G}_d = \sum_{k=0}^{N-1} \sum_{i=1}^{s} & h \left[ b_i L(Q_k^i, V_k^i)  + \left\langle p_k^i, (Q_k^i - q_k)/h -   \sum_{j=1}^s a_{ij} V_k^j \right\rangle  \right.  \\
&+ \left. \left\langle p_{k+1}, (q_{k+1} - q_k)/h -   \sum_{j=1}^s b_{j} V_k^j \right\rangle +  b_i \left\langle \Lambda_k^i, g(Q_k^i) \right\rangle \right]  \text{.}
\end{align*}
\end{definition}

The following condition on the coefficients of the s-stage RK method will be important in obtaining a well-defined discrete update on phase space that also respects the constraints \citep{Ja1996}.

\begin{Condition} \label{cond:RKcoefficients} Consider an s-stage RK method with given $b$-vector and $a$-matrix.  Assume $b_k \ne 0$ and set $\hat{a}_{kj} =  b_j - a_{j k} b_j/b_k$.  The coefficients of the s-stage RK method satisfy:
\begin{align*}
& a_{1i}=0,~~~a_{si} = b_i,~~~b_i \ne 0,~~i=1,...,s \\
& \left(\sum_{k=1}^s a_{ik} \hat{a}_{kj} \right)_{i,j=2}^s~~\text{invertible}\text{.}
\end{align*}
\end{Condition}

Under this condition on the coefficients of the s-stage RK method, one can prove the following theorem.

\begin{theorem} \label{thm:cvprk}
Given an s-stage RK method with $b$-vector and $a$-matrix  that satisfy condition~\ref{cond:RKcoefficients}, a Lagrangian system with smooth Lagrangian $L: TQ \to \mathbb{R}$ such that $\partial^2 L/ \partial v^2$ is invertible, and smooth holonomic constraint function $g: Q \to \mathbb{R}^k$.     A discrete curve $c_d \in \mathcal{C}_d(q_1, q_2)$ satisfies the following VPRK method:
\begin{equation}  \label{eq:cvprk}  
\begin{cases}
\begin{array}{rll}
Q_k^i &=  & q_k + h \sum_{j=1}^s a_{ij} V_k^j,  \\
q_{k+1} &=  &q_k + h \sum_{j=1}^s b_j V_k^j,   \\
P_{k}^i &= & p_k + h \sum_{j=1}^s \left( b_j  - \frac{b_j a_{ji}}{b_i} \right) \left( \frac{\partial L}{\partial q}(Q_k^j, V_k^j) + \frac{\partial g}{\partial q}(Q_k^j)^* \cdot \Lambda_k^j \right),  \\
p_{k+1} &= &p_k + h \sum_{j=1}^s b_j \left( \frac{\partial L}{\partial q}(Q_k^j, V_k^j)  + \frac{\partial g}{\partial q}(Q_k^j)^* \cdot \Lambda_k^j \right),   \\
P_k^i &= & \frac{\partial L}{\partial v} (Q_k^i, V_k^i), \\
g(Q^i_k)  &=& 0   \text{.}
\end{array}
\end{cases}
\end{equation}
for $i=1,\cdots,s$ and $k=1,\cdots,N-1$, if and only if it is a critical point of the function $\mathfrak{G}_d:  \mathcal{C}_d \to \mathbb{R}$, that is, $\mathbf{d} \mathfrak{G}_d (c_d) = 0$.  Moreover, there exist $\{ \Lambda_k^i \}_{i=1}^s$ such that the discrete flow map defined by the above scheme, $F_h: T^*S \to T^*S$, preserves the canonical symplectic form on $T^*S$.  
\end{theorem}

\begin{proof}
Under the assumptions of the theorem, the existence of a numerical solution of (\ref{eq:cvprk}) and a discrete flow map $F_h:T^*S \to T^*S$ is guaranteed.  That is, given $(q_k, p_k) \in T^*S$ one can solve (\ref{eq:cvprk}) for $(q_{k+1}, p_{k+1}) \in T^*S$.  In particular, the condition that $\partial^2 L/ \partial v^2$ is invertible, implies one can eliminate $\{ V_k^i \}_{i=1}^s$ using the Legendre transform of $L$.  The condition~\ref{cond:RKcoefficients} implies that one can determine $\{ \Lambda_k^i \}_{i=1}^s$ so that $(q_{k+1}, p_{k+1}) \in T^*S$ \citep{Ja1996}.

The differential of $\mathfrak{G}_d(c_d)$  in the direction $z = (\{ \delta q_k, \delta p_k \}, \{ \delta Q_k^i, \delta V_k^i, p_k^i \}_{i=1}^s, \{ \delta \Lambda_k^i \}_{i=1}^s)$ is given by:
\begin{align*}
\mathbf{d} &\mathfrak{G}_d \cdot z =
\sum_{k=0}^{N-1} \sum_{i=1}^{s} h b_i \left[ \frac{ \partial L}{\partial q}(Q_k^i, V_k^i) \cdot \delta Q_k^i + \frac{ \partial L}{\partial v}(Q_k^i, V_k^i) \cdot \delta V_k^i  \right] \\
&+ h \left[ \left\langle p_k^i, (\delta Q_k^i - \delta q_k)/h -\sum_{j=1}^s a_{ij} \delta V_k^j \right\rangle  + \left\langle p_{k+1}, (\delta q_{k+1} - \delta q_k)/h -   \sum_{j=1}^s b_{j} \delta V_k^j \right\rangle \right]  \\
&+ h \left[ \left\langle \delta p_k^i, (Q_k^i - q_k)/h -\sum_{j=1}^s a_{ij} V_k^j \right\rangle  + \left\langle \delta p_{k+1}, (q_{k+1} - q_k)/h -   \sum_{j=1}^s b_{j} V_k^j \right\rangle \right]  \\
&+ h b_i \left[ \left\langle \delta \Lambda_k^i, g(Q_k^i) \right\rangle +  \left\langle \Lambda_k^i, \frac{\partial g}{\partial q}(Q_k^i) \cdot \delta Q_k^i \right\rangle \right] \text{.}
\end{align*}
Collecting terms with the same variations and summation by parts using the boundary conditions $\delta q_0 = \delta q_N = 0$ gives,
\begin{align*}
\mathbf{d} &\mathfrak{G}_d \cdot z =
\sum_{k=1}^{N-1}   \sum_{i=1}^s  \left( h b_i \frac{\partial g}{\partial q}(Q_k^i)^* \Lambda_k^i  + h b_i \frac{ \partial L}{\partial q}(Q_k^i, V_k^i) + p_k^i \right) \cdot \delta Q_k^i \\
&+  \left( -p_{k+1} + p_k - \sum_{i=1}^s p_k^i  \right) \cdot   \delta q_k +h \left( b_i \frac{ \partial L}{\partial v}(Q_k^i, V_k^i) - \sum_{j=1}^s a_{ji} p_k^j  - b_i p_{k+1}  \right) \cdot \delta V_k^i \\
&+ h  \left\langle \delta p_k^i, (Q_k^i - q_k)/h -\sum_{j=1}^s a_{ij} V_k^j \right\rangle  + h \left\langle \delta p_{k+1}, (q_{k+1} - q_k)/h -   \sum_{j=1}^s b_{j} V_k^j \right\rangle \\
&+ h b_i \left\langle \delta \Lambda_k^i, g(Q_k^i) \right\rangle  \text{.}
\end{align*}
Since $\mathbf{d} \mathfrak{G}_d(c_d) =0$ if and only if $\mathbf{d} \mathfrak{G}_d \cdot z = 0$ for all $z \in T_{c_d} \mathcal{C}_d$, one arrives at the desired equations with the elimination of $p_k^i$ and the introduction of the internal stage variables $P_k^i = \partial L/\partial v (Q_k^i, V_k^i) $ for $i=1,\cdots,s$.   Conversely, if $c_d$ satisfies (\ref{eq:cvprk}) then $\mathbf{d} \mathfrak{G}_d(c_d) =0$.

We will employ the variational proof of symplecticity to prove that $F_h$ is symplectic.  This proof is standard, however, we provide it here to emphasize that the symplectic form that is exactly preserved by the method is the canonical one on $T^*S$.

Consider the subset of $\mathcal{C}_d$ given by solutions of (\ref{eq:cvprk}).  Let $\hat{\mathfrak{G}}_d$ denote the restriction of $\mathfrak{G}_d$ to this space. Since each of these solutions is determined by an initial point on $T^*S$, one can identify this space with $T^*S$, and hence, $\hat{\mathfrak{G}}_d : T^*S \to \mathbb{R}$.   Since $\hat{\mathfrak{G}}_d$ is restricted to solution space,
\begin{align*}
\mathbf{d} &\hat{\mathfrak{G}}_d(q_0, p_0) \cdot (\delta q_0, \delta p_0) = \left\langle p_N, \delta q_N \right\rangle - \left\langle p_0, \delta q_0 \right\rangle
  \text{.}
\end{align*}
Observe that these boundary terms are the canonical one-forms on $T^*S$ evaluated at $k=0$ and $k=N$.  Preservation of the symplectic form follows from $\mathbf{d}^2 \hat{\mathfrak{G}}_d = 0$.
\end{proof}

\paragraph{Variational RATTLE for Constrained Systems}

The variational RATTLE integrator is the Lagrangian analog of the RATTLE algorithm originally proposed as a constrained version of Verlet in \citep{RyCiBe1977}.  It was shown to be symplectic in \citep{LeSk1994}, and was extended to general constrained Hamiltonian systems by \citep{Ja1996}.   It is defined by the following two-stage RK discretization of the kinematic constraint (implicit trapezoidal rule),
\begin{center}
\begin{tabular}{c|cc}
 &  $0$  & $0$ \\
  & $1/2$    & $1/2$ \\
\hline  
          & $1/2$    & $1/2$ 
\end{tabular}
\end{center}
Given $h$ and $(q_k, p_k)$, the method determines $(q_{k+1}, p_{k+1})$ and two Lagrange multipliers $\Lambda_k^{1}, \Lambda_k^{2}$ by solving the following system of equations, 
\begin{equation} \label{eq:rattle}
\begin{cases}
\begin{array}{ccl}
q_{k+1} &=& q_k + \frac{h}{2}  \left( V_k^1 +  V_k^2 \right)  \text{,} \\
P_k^1 &=& p_k + \frac{h}{2} \left( \frac{\partial L}{\partial q}(q_k, V_k^1) + \frac{\partial g}{\partial q}(q_k)^* \Lambda^1_k \right) \text{,} \\
p_{k+1} &=& P_k^1 + \frac{h}{2} \left( \frac{\partial L}{\partial q}(q_{k+1}, V_k^2) + \frac{\partial g}{\partial q}(q_{k+1})^* \Lambda^2_k \right)  \\
0 &=& g(q_{k+1})  \text{,}  \\
0 &=& \frac{\partial g}{\partial q}(q_{k+1}) \cdot v_{k+1} \text{,} \\
p_{k+1} &=& \frac{\partial L}{\partial v}(q_{k+1}, v_{k+1})   \text{,}\\
P_k^1 &=& \frac{\partial L}{\partial v}(q_k, V_k^1)  = \frac{\partial L}{\partial v}(q_{k+1}, V_k^2) \text{.} 
\end{array}
\end{cases}
\end{equation}   
By theorem~\ref{thm:cvprk}, variational RATTLE defines a symplectic scheme.  Moreover, as is well-known it is second-order accurate.

\paragraph{Variational Euler for Constrained Systems}

One can relax the conditions assumed in theorem \ref{thm:cvprk} on the coefficients of the s-stage RK method \citep{Ha2003}.   For instance, if $b_s=0$ in the s-stage RK method one can still obtain a well-defined variational integrator.   Consider as an example the following two-stage RK method,
\begin{center}
\begin{tabular}{c|cc}
 &  $0$  & $0$ \\
  & $1$  & $0$ \\
\hline  
          & $1$    & $0$ 
\end{tabular}
\end{center}
The corresponding variational integrator is given by:
\begin{equation} \label{eq:variationaleulerA}
\begin{cases}
\begin{array}{ccl}
q_{k+1} &=& q_k + h  \hat{v}_{k+1}  \text{,} \\
\hat{p}_{k+1} &=& p_k + h \left( \frac{\partial L}{\partial q}(q_k, \hat{v}_{k+1}) + \frac{\partial g}{\partial q}(q_k)^* \Lambda^1_k \right)  \\
0 &=& g(q_{k+1})  \text{,}  \\
\hat{p}_{k+1} &=& \frac{\partial L}{\partial v}(q_k, \hat{v}_{k+1})   \text{.}
\end{array}
\end{cases}
\end{equation}  
However, the corresponding discrete flow does not satisfy the ``hidden'' velocity constraint, and hence, does not define a map from $T^*S$ to $T^*S$.  To satisfy the hidden constraint a projection step $(q_{k+1}, \hat{p}_{k+1}) \mapsto (q_{k+1}, p_{k+1})$ is taken:
\begin{equation} \label{eq:projectionstep}
\begin{cases}
\begin{array}{ccl}
p_{k+1} &= & \hat{p}_{k+1} + h \frac{\partial g}{\partial q}(q_{k+1})^* \Lambda^2_k \text{,} \\
0 &=& \frac{\partial g}{\partial q}(q_{k+1}) \cdot v_{k+1} \text{,} \\
p_{k+1} &=& \frac{\partial L}{\partial v}(q_{k+1}, v_{k+1})   \text{.}
\end{array}
\end{cases}
\end{equation}
One can check that this projection step defines a symplectic map.   Thus, the composite map $(q_k, p_k) \mapsto (q_{k+1}, \hat{p}_{k+1})  \mapsto (q_{k+1}, p_{k+1})$ is symplectic.  It is the Lagrangian version of the constrained symplectic Euler method \citep{HaLuWa2006}.

Another example relaxing the assumptions in theorem \ref{thm:cvprk} is given by regarding the 2-stage RK method as being single-stage implicit Euler,
\begin{center}
\begin{tabular}{c|c}
  & $1$   \\
\hline  
          & $1$    
\end{tabular}
\end{center}
The corresponding variational integrator is:
\begin{equation}  \label{eq:variationaleulerB}
\begin{cases}
\begin{array}{ccl}
q_{k+1}& =& q_k + h  \hat{v}_{k+1}  \text{,} \\
\hat{p}_{k+1}& =& p_k + h \left( \frac{\partial L}{\partial q}(q_k, v_k) + \frac{\partial g}{\partial q}(q_k)^* \Lambda^1_k \right)  \\
0 &=& g(q_{k+1})  \text{,}  \\
\hat{p}_{k+1} &=& \frac{\partial L}{\partial v}(q_{k+1}, \hat{v}_{k+1})   \text{.} 
\end{array}
\end{cases}
\end{equation}
We conclude with a theorem summarizing the structure-preserving properties of these constrained variational Euler methods.

\begin{theorem}
The composition of one step of (\ref{eq:variationaleulerA}) or (\ref{eq:variationaleulerB}) and the projection step (\ref{eq:projectionstep}), defines a symplectic integrator on $T^*S$.  Moreover, these integrators are first-order accurate.
\end{theorem}

\section{Constrained, Stochastic Mechanical Systems}  \label{sec:constraintedstochastichpmechanics}

\paragraph{Constrained, Stochastic HP Principle}  
The setting in this section is an n-manifold $Q$ and a {\em stochastic mechanical system} with smooth, holonomic constraint function, $g: Q \to \mathbb{R}^k$, $k<n$, that has a regular value $0 \in \mathbb{R}^k$.  Its configuration space is given by the {\em constraint submanifold}: $S = g^{-1}(0)$.   In this section we will introduce Lagrange multipliers to prove an equivalence between a stochastic variational principle on $S$ and a constrained, stochastic variational principle on $Q$.

We will adopt a HP viewpoint to describe this mechanical system with random perturbations and will refer to this perturbed system as a {\bfi stochastic mechanical system}.  The HP principle unifies the Hamiltonian and Lagrangian descriptions of a mechanical system  \citep{YoMa2006a, YoMa2006b, Bo2007, BoMa2007}.    Roughly speaking, in the stochastic context it states the following critical point condition on $PS$ (cf.~definition~\ref{defn:hppathspace}), 
\[
\delta \int _a^b \left[ L(q, v) dt + \sum_{i=1}^m \gamma_i (q) \circ d W_i + \left\langle p, \frac{dq}{dt} - v \right\rangle dt \right]  = 0,
\]  
where $(q(t),v(t),p(t)) \in PS$ are varied arbitrarily and independently with endpoint conditions $q (a) $ and $q (b)$ fixed, and $\gamma_i: Q \to \mathbb{R}$ and $W_i$ for $i=1,...,m$ are {\bfi stochastic potentials} and Wiener processes respectively. This principle builds in a Legendre transform, stochastic Hamilton's equations and stochastic Euler--Lagrange equations.    The action integral in the above principle, consists of two Lebesgue integrals with respect to the Lebesgue measure $dt$ and $m$ Stratonovich stochastic integrals.  This action is random, i.e., for every sample point $\omega \in \Omega$ one will obtain a different, time-dependent Lagrangian system.   However, each system possesses a variational structure as made precise in \citep{BoOw2007a}.    The following definitions will be useful to state the constrained, variational principle of HP for mechanical systems with holonomic constraints.

\begin{definition}
Let $i: S \to Q$ be the inclusion mapping. Fixing the interval $[a,b]$, define {\bfi constrained HP path space} as
 \begin{align*}
\mathcal{C}_c^{HP}([a,b],q_1,q_2) =& \{ (q,v,p): [a,b] \to PS ~|~ \\
& q \in C^1([a,b], S),~(v,p) \in C^0([a,b]),~q(a)=q_1,~q(b)=q_2 \}  \text{,}
\end{align*}
and the {\bfi unconstrained HP path space} as
 \begin{align*}
\mathcal{C}^{HP}([a,b],q_1,q_2) =& \{ (q,v,p): [a,b] \to PQ ~|~ \\
& q \in C^1([a,b],Q),~(v,p) \in C^0([a,b]),~q(a)=i(q_1),~q(b)=i(q_2) \}  \text{.}
\end{align*}
\end{definition}

Let $(\Omega, \mathcal{F}, P)$ be a probability space, $\mathcal{F}_t$, $t \in [a,b]$, be a nondecreasing family of $\sigma$-subalgebras of $\mathcal{F}$, $\omega \in \Omega$, and $(W_i(t, \omega), \mathcal{F}_t)$, $i=1,...,m$, be independent, real-valued Wiener processes.   In terms of these Wiener processes, we define the following.

\begin{definition}
Set $L^S = \left. L \right|_{TS}$ and $\gamma_j^S = \left. \gamma_j \right|_{S}$ for $j=1,...,m$. Moreover define the {\bfi unconstrained action} $\mathfrak{G}: \Omega \times \mathcal{C}([a,b],q_1,q_2) \to \mathbb{R}$  as
\[
\mathfrak{G} =  \int _a^b \left[ L(q, v) dt + \sum_{j=1}^m \gamma_j (q) \circ d W_j + \left\langle p, \frac{dq}{dt} - v \right\rangle dt \right] \text{.}
\]
and the constrained action as $\mathfrak{G}_c = \left. \mathfrak{G}  \right|_{\mathcal{C}_c^{HP}}$.
\end{definition}

The unconstrained HP path space is a smooth infinite-dimensional manifold.  One can show that its tangent space at  $c = (q,v,p) \in \mathcal{C}^{HP}( [a,b], q_1, q_2)$ consist of maps $w = (q,v,p, \delta q, \delta v, \delta p) \in C^0([a,b],T(PQ))$  such that  $\delta q(a) = \delta q(b) = 0$ and $q, \delta q$ are of class $C^1$.

As opposed to using generalized coordinates on $PS$, we wish to describe the mechanical system using constrained coordinates on $PQ$ and introduce Lagrange multipliers to enforce the constraint.  However, because of the stochastic component of the action, the standard Lagrange multiplier theorem will not apply directly and one cannot introduce Lagrange multipliers in the standard way.  Instead, we will introduce the Lagrange multiplier using the following.

\begin{definition}
Given $f_1 \in C^0([a,b], T^*Q)$ and $f_2 \in C^1([a,b], TQ)$, define
\[
\int_a^t \left\langle df_1,  f_2 \right\rangle :=  
\left. \left\langle f_1, f_2 \right\rangle \right|_{a}^t - \int_a^t \left\langle f_1,  \frac{d f_2}{dt} \right\rangle dt,~~ t \in [a,b]   \text{.}
\]
\label{defn:dfintegral}
\end{definition}

Differentiability of the flow map on $T^*S$ will be defined in the mean-squared sense.  In the following we define mean-squared derivatives on a Banach space $E$ with the understanding that this notion can be extended to any manifold using a local representative of the flow map.

\begin{definition}[Mean-squared Derivative]
The mean squared norm of $f: E \times \Omega \to E$ is given by:
\[
\| f (x, \omega) \| = \left( \E \left( |f(x, \omega)|^2 \right) \right)^{1/2}
\]
Using this norm one can define the derivative of $f$ in the standard way, i.e., $f$ is mean squared differentiable at $a \in E$ if there is a bounded linear map $Df(a, \omega): E \to E$ that satisfies,
\[
\lim_{\delta \to 0} \frac{ \left\| f(a+\delta, \omega) - f(a, \omega) - D f(a, \omega) \cdot \delta \right\| }{ \| \delta \| } \to 0  \text{.}
\]
\label{def:meansquared}
\end{definition}

With the above definitions we prove the following.

\begin{theorem}[{\bf Constrained, Stochastic HP Principle}]  Given a stochastic mechanical system with Lagrangian $L: TQ \to \mathbb{R}$ such that $\partial^2 L/ \partial v^2$ is invertible, stochastic potentials $\gamma_j: Q \to \mathbb{R}$ for $j=1,...,m$, and holonomic constraint $g: Q \to \mathbb{R}^k$ with $S = g^{-1}(0)$. Then the following are equivalent:
\begin{description}
\item[(i)] $z=(q, v, p) \in \mathcal{C}_c^{HP}([a,b],q_1,q_2)$ extremizes $\mathfrak{G}_c$.
\item[(ii)] $z=(q, v, p) \in \mathcal{C}_c^{HP}([a,b],q_1,q_2)$ satisfies stochastic HP equations
\begin{equation} \label{eq:stochastichp}
\begin{cases}
\begin{array}{ccl}
dq &=& v dt   \text{,}  \\
d p &= & \frac{\partial L^S}{\partial q}(q,v) dt + \sum_{j=1}^m \frac{\partial \gamma_j^S }{\partial q}(q) \circ d W_j \text{,}  \\
p &=& \frac{\partial L^S}{\partial v}(q,v) \text{.}
\end{array}
\end{cases}
\end{equation}
\item[(iii)] There exists $\lambda \in C^0([a,b],\mathbb{R}^k)$ such that $z=(q, v, p) \in \mathcal{C}^{HP}([a,b],i(q_1),i(q_2))$ and $\lambda$ extremize the augmented action $\bar{\mathfrak{G}}(z,\lambda) =  \mathfrak{G}(z) + \left\langle \lambda, \Phi(z) \right\rangle$  where $\Phi(z)(t) = dg/dt(q(t))$ and $\left\langle \lambda, \Phi(z) \right\rangle = \int_a^b \left\langle \lambda(s), \Phi(z)(s) \right\rangle ds$.   
\item[(iv)] There exists $\lambda \in C^0([a,b],\mathbb{R}^k)$ such that $z=(q, v, p) \in \mathcal{C}^{HP}([a,b],i(q_1),i(q_2))$ and $\lambda$ satisfy the constrained, stochastic HP equations
\begin{equation} \label{eq:cstochastichp}
\begin{cases}
\begin{array}{ccl}
dq &=& v dt \text{,}  \\
d p &=&  \frac{\partial L}{\partial q}(q,v) dt + \sum_{j=1}^m \frac{\partial \gamma_j}{\partial q}(q) \circ d W_j + \frac{\partial g}{\partial q}(q)^* \cdot d \lambda \text{,}  \\
p &=& \frac{\partial L}{\partial v}(q,v) \text{,}   \\
g(q) &=& 0  \text{.}  
\end{array}
\end{cases}
\end{equation}
\end{description}
Moreover, the flows of (\ref{eq:stochastichp}) and (\ref{eq:cstochastichp}) are mean-square symplectic.
\label{thm:stochasticconstrainedhpprinciple}
\end{theorem}

\begin{Remark}
The constrained, stochastic HP equations should be thought of in integral form.  In particular, the Lagrange multipler term in (\ref{eq:cstochastichp}) by definition~\ref{defn:dfintegral} satisfies:
\[
\int_a^t \left\langle \frac{\partial g}{\partial q}(q)^* \cdot d \lambda, w \right\rangle   =
\left.   \left\langle \lambda, \frac{\partial g}{\partial q}(q) w \right\rangle  \right|_{a}^t  - \int_a^t  \left\langle  \lambda, \frac{d}{dt}  \frac{\partial g}{\partial q}(q) w \right\rangle dt \text{,}
\]
for any $w \in \mathbb{R}^k$.
\end{Remark}

\begin{proof}
The stochastic HP principle states that (i) and (ii) are equivalent \citep{BoOw2007a}.  Assume that (iii) is true.  Then for all $v_z = (v_q, v_v, v_p) \in T_z \mathcal{C}^{HP}$,
\begin{equation}
\mathbf{d} \bar{\mathfrak{G}} (z, \lambda) \cdot (v_z, v_{\lambda} ) = \mathbf{d} \mathfrak{G}(z) \cdot v_z +
\left\langle v_{\lambda}, \Phi(z) \right\rangle + \left\langle \lambda, \mathbf{d} \Phi(z) \cdot v_{z} \right\rangle = 0 
\label{eq:cpbarG}
\end{equation}
Since $\Phi(z) = 0$ for $ z \in \mathcal{C}_c^{HP}$ and $T_{q(t)} g(q(t)) \cdot v_q(t) = 0$ for $v_q(t) \in T_q S$, 
\[
\mathbf{d} \bar{\mathfrak{G}} (z, \lambda) \cdot (v_z, v_{\lambda} ) =\mathbf{d} \mathfrak{G}(z) \cdot v_z = 0,~~\forall~~v_z \in T_z \mathcal{C}_c^{HP} 
\]
which implies (i) and (ii).

Moreover, expanding (\ref{eq:cpbarG}) and setting $(v_q, v_v, v_p) = (\delta q, \delta v, \delta p)$ and $v_{\lambda} = \delta \lambda$ yields,
\begin{align*}
\mathbf{d} & \bar{\mathfrak{G}}  (z, \lambda)  \cdot  (\delta q, \delta v, \delta p, \delta \lambda) = \int_a^b \left[ \left\langle \frac{\partial L}{\partial q} , \delta q \right\rangle ds + \sum_{i=1}^m \left\langle \frac{\partial \gamma_i}{\partial q},  \delta q \right\rangle \circ  d W_i +  \left\langle \frac{\partial L}{\partial v} , \delta v \right\rangle ds \right. \\
& \left.  + \left\langle \delta p,  \frac{d q}{dt} - v  \right\rangle ds + \left\langle p, \delta \frac{d q}{dt} - \delta v  \right\rangle ds + \left\langle \delta \lambda, \frac{d}{dt} g(q) \right\rangle ds
+ \left\langle \lambda, \frac{d}{dt} \left( \frac{\partial g}{\partial q} \cdot \delta q \right) \right\rangle ds \right] = 0\text{.}
\end{align*}
Consider the terms involving $\delta p$, $\delta v$ and $\delta \lambda$.  Since these variations are arbitrary, the following hold:\footnote{This follows from the basic lemma that if $f, g \in C^0([a,b], \mathbb{R})$ and $g$ is arbitrary then $\int_a^b f(t) g(t) dt = 0 \iff f(t) = 0~~\forall~~t \in[a,b]$.}
\[
\frac{dq}{dt} = v,~~~\frac{\partial L}{\partial v}(q,v) = p,~~~\frac{d}{dt} g(q) = 0.
\]
Since $g(q(a)) = g(i(q_1)) = 0$, $d/dt ( g(q) ) = 0$ implies that $g(q(t)) = 0$ for all $t \in [a,b]$.

Collecting the variations with respect to $\delta q$ in the differential gives,
\[
\int_a^b \left[ \frac{\partial L}{\partial q}  \cdot \delta q ds + \left\langle p, \delta \frac{d q}{dt} \right\rangle ds + 
  \sum_{j=1}^m \frac{\partial \gamma_j}{\partial q} \cdot \delta q \circ  d W_j  + \left\langle \lambda, \frac{d}{dt} \left( \frac{\partial g}{\partial q} \cdot \delta q \right) \right\rangle ds \right] = 0 
  \text{.}
\]
Using definition~\ref{defn:dfintegral}, we introduce the following function $I: \mathcal{C}^{HP}(q_1, q_2, [a,b]) \times C([a,b], \mathbb{R}^k) \times C^1([a,b], TQ)  \to \mathbb{R}$,
\[
I(q,v,p,\lambda,f) = \int_a^b \left[ \left( \frac{\partial L}{\partial q} ds  + \sum_{j=1}^m \frac{\partial \gamma_j}{\partial q} \circ d W_j  - dp - \frac{\partial g}{\partial q}^* \cdot d \lambda \right) \cdot f  \right]  \text{.}
\]
In the following it is shown that if $I(q,v,p, \lambda, f) = 0$ for arbitrary $f$ of class $C^1$ then $c=(q,v,p,\lambda)$ satisfies (\ref{eq:cstochastichp}).

Let $\{ U_{\alpha}, g_{\alpha} \}$ be a partition of unity on $PQ \times \mathbb{R}^k$.  Expand $I$ in terms of this partition of unity,
\[
I = \sum_{\alpha} \int_a^b \left[ g_{\alpha}(q,v,p,\lambda) \left( \frac{\partial L}{\partial q} dt  + \sum_{j=1}^m \frac{\partial \gamma_j}{\partial q} \circ  d W_j  - dp - \frac{\partial g}{\partial q}^* \cdot d \lambda \right) \cdot f \right]  \text{.}
\]
Since the curves $z=(q,v,p)$ and $\lambda$ are compactly supported, only a finite number of the $g_{\alpha}$ are nonzero.  For each $g_{\alpha}$ nonzero, the terms in the integral can be expressed in local coordinates.     Observe that since $dq = v dt$, the Stratonovish-Ito conversion formula implies that, 
\[
 \int_a^b  g_{\alpha} \frac{\partial \gamma_j}{\partial q} \cdot \delta q  \circ d W_j  =
  \int_a^b  g_{\alpha} \frac{\partial \gamma_j}{\partial q} \cdot \delta q  d W_j
\]
for $j=1,...,m$.

We will select $f$ to single out the $\beta$th-component of the covector field in $I$.  Introduce the following function $h: \mathbb{R} \to \mathbb{R}$ for this purpose:
\[
h(t) =  2 \frac{t}{\epsilon} - \frac{t^2}{\epsilon^2} \text{.}
\] 
Observe that $h(0)=0$, $h(\epsilon) = 1$, and $h^{\prime}(\epsilon)=0$.  Let $\{ e_{\beta} \}$ be a basis for the model space of $Q$.  Now fix $\beta$, and define $f_{\epsilon} \in C^1([a,b], TQ)$ in local coordinates as:
\[
f_{\epsilon}(s)  = 
\begin{cases}
h(s-a) e_{\beta} &  \text{if}~ ~a \le s \le a + \epsilon \text{,} \\
e_{\beta}  &  \text{if}~~ a + \epsilon < s < t - \epsilon \text{,} \\
h(t-s) e_{\beta} &  \text{if}~~ t - \epsilon  \le s \le t  \text{,} \\
 0 &  \text{if}~~  t  < s \le b \text{.}
\end{cases}
\]
Introduce the following label to simplify subsequent calculations,
\[
A(s)  = \left( \frac{\partial L}{\partial q}(q(s), v(s)) ds  + \sum_{i=1}^m \frac{\partial \gamma_i}{\partial q}(q(s)) ~d W_i(s) - dp(s) - \frac{\partial g}{\partial q}(q(s))^* \cdot d \lambda(s) \right) \cdot e_{\beta}  \text{.}
\]
In terms of $A(s)$,  one can write
\[
I(q,v,p,\lambda, f_{\epsilon}) = \sum_{\alpha} \left[ \int_{a}^{a+\epsilon} h(s-a) g_{\alpha}(s) A(s) 
+ \int_{a+\epsilon}^{t-\epsilon} g_{\alpha}(s) A(s) 
+ \int_{t-\epsilon}^t  h(t-s) g_{\alpha}(s) A(s) 
\right]  \text{.}
\]

We will show in the mean squared norm,
\begin{equation}
\lim_{\epsilon \to 0}  I(q,v,p,\lambda,f_{\epsilon})  = \sum_{\alpha} \int_{a}^{t} g_{\alpha} A(s) =: I^*  \text{.}
\label{eq:msIconvergence}
\end{equation}
Using this result and the Borel-Cantelli lemma, one can deduce there exists $\{ \epsilon_n \}$ that converges to $0$  such that $I(q,v,p,\lambda,f_{\epsilon_n})$ a.s.~converges to $I^*$.  It follows that $I^* = 0$ almost surely.  


We proceed to prove (\ref{eq:msIconvergence}).  Since $(a+b)^2 \le 2 a^2 + 2 b^2$,
\begin{align*}
 & \left\|  \sum_{\alpha} \int_{a}^{t} g_{\alpha} A(s) -  I(q,v,p,\lambda,f_{\epsilon})  \right\|^2  \\
& = \left\|  \sum_{\alpha}
\int_{a}^{a+\epsilon} ( 1- h(s-a)) g_{\alpha} A(s) 
+ \int_{t-\epsilon}^t (1- h(t-s)) g_{\alpha} A(s) \right\|^2  \\
& \le   2 \left\| \sum_{\alpha} \int_{a}^{a+\epsilon} (1-h(s-a)) g_{\alpha} A(s) \right\|^2 
+ 2 \left\| \sum_{\alpha} \int_{t-\epsilon}^t  (1-h(t-s))  g_{\alpha} A(s) \right\|^2
\text{.}
\end{align*}
We will only show how to bound the first term since bounding the second term is very similar.  By continuity of $(q,v,p,\lambda)$, one can pick $\epsilon>0$ small enough so that the support of $(q,v,p,\lambda)$ lies in a single chart.  Therefore,
\begin{align*}
 & \left\|  \sum_{\alpha} \int_{a}^{a+\epsilon} (1-h(s-a)) g_{\alpha} A(s) \right\|^2  \\
& = \left\| \int_{a}^{a+\epsilon} (1-h(s-a))  
\left( \frac{\partial L}{\partial q} ds  + \sum_{j=1}^m \frac{\partial \gamma_j}{\partial q} ~d W_j - dp - \frac{\partial g}{\partial q}^* \cdot d \lambda \right) \cdot e_{\beta}  \right\|^2  \\
&\le 4 \left\| \int_{a}^{a+\epsilon}  (1-h(s-a))  \frac{\partial L}{\partial q^{\beta}} ds \right\|^2   +  4 \sum_{j=1}^m \left\| \int_{a}^{a+\epsilon} (1-h(s-a)) \frac{\partial \gamma_j}{\partial q^{\beta}} ~d W_j \right\|^2  \\
 &+ 4 \left\|   \int_{a}^{a+\epsilon} (1-h(s-a)) dp \cdot e_{\beta} \right\|^2 +
  4 \left\| \int_{a}^{a+\epsilon} (1-h(s-a)) \left\langle \frac{\partial g}{\partial q}^* \cdot d \lambda, e_{\beta} \right\rangle \right\|^2 \text{.}
\end{align*}
Since $\frac{\partial L}{\partial q^j}$ is continuous on $s \in[a,a+\epsilon]$, the first term can be bounded,
\[
\left\| \int_{a}^{a+\epsilon} (1-h(s-a))   \frac{\partial L}{\partial q^{\beta}} ds \right\|^2 \le
\frac{M^2  \epsilon^2}{9} \text{.}
\]
Similarly, by the Ito isometry and since $\frac{\partial \gamma_i}{\partial q^{\beta}}$ is continuous on $s \in[a, a +\epsilon]$, the second $m$ terms can similarly be bounded, e.g., the jth Stratonovich integral can be bounded as follows,
\[
 \left\| \int_{a}^{a+\epsilon} (1- h(s-a)) \frac{\partial \gamma_j}{\partial q^{\beta}} ~d W_j\right\|^2  =
  E \left( \int_{a}^{a+\epsilon} \left| (1-h(s-a)) \frac{\partial \gamma_j}{\partial q^{\beta}} \right|^2 ds \right) \le \frac{M^2 \epsilon}{5} \text{.}
\]
Using (\ref{defn:dfintegral}) and the integral mean value theorem, there exists a real constant $c_1 \in [0, 1]$ such that:
\begin{align*}
\left\|   \int_{a}^{a+\epsilon} (1-h(s-a)) dp \cdot e_{\beta} \right\|^2 & = 
 \left\| \left. (1-h(s-a)) p_{\beta} (s) \right|_{a}^{a+\epsilon} + \int_a^{a+\epsilon} p_{\beta}(s) h^{\prime}(s-a) ds \right\|^2  \\
 & =   \left\| - p_{\beta}(a)  + p_{\beta}(a + c_1 \epsilon)  \right\|^2  \text{.}
\end{align*}
Similarly, there exist constants $c_2, c_3 \in [0,1]$ such that
\begin{align*}
& \left\|   \int_{a}^{a+\epsilon} (1-h(s-a))  \left\langle \frac{\partial g}{\partial q}^* \cdot d \lambda, e_{\beta} \right\rangle \right\|^2 = \\
 =&  \left\| \left. (1-h(s-a)) \left\langle \frac{\partial g}{\partial q}^* \cdot \lambda, e_{\beta} \right\rangle \right|_{a}^{a+\epsilon} - \int_a^{a+\epsilon} ( 1- h(s-a) ) \lambda_{\beta} \frac{d}{ds} \frac{\partial g}{\partial q^{\beta}}(q(s)) ds \right. \\
& \left. + \int_a^{a+\epsilon} \left\langle   \lambda,  \frac{\partial g}{\partial q} \cdot e_{\beta} \right\rangle h^{\prime}(s-a) ds\right\|^2  \\
 =&   \left\| - \left\langle \frac{\partial g}{\partial q}(q(a))^* \cdot \lambda(a), e_{\beta} \right\rangle  
 - \frac{\epsilon}{3} \lambda_{\beta}(a+c_2 \epsilon)  \frac{\partial g}{\partial q^{\beta}}(q(a + c_2 \epsilon))  \right. \\
& \left. + \left\langle   \frac{\partial g}{\partial q}(q(a+c_3 \epsilon))^* \cdot \lambda(a+c_3 \epsilon),   e_{\beta}  \right\rangle
  \right\|^2  \text{.}
\end{align*} 
Since $p_{\beta}$ and $\lambda$ are of class $C^0$ and $g$ is smooth, as $\epsilon \to 0$ these terms vanish.  Since $\beta$ is arbitrary we have proved (\ref{eq:msIconvergence}).  Therefore, almost surely: if $c=(z, \lambda)$ is a critical point of $\mathfrak{G}$ then $c$ satisfies (\ref{eq:cstochastichp}).

On the other hand, assume (iv) is true.  Then almost surely: if $c=(z, \lambda)$ satisfies the constrained, stochastic HP equations, then it is a critical point of $\bar{\mathfrak{G}}$.  This direction is easy to confirm, since as a solution to the constrained, stochastic HP equations $c$ is a measureable diffusion process.   In fact, this direction is similar to the one Bismut originally established, namely that the solution of stochastic Hamilton's equations extremize an action function; albeit a different stochastic action is used in this proof \citep{Bi1981}.

Assume that (i) is true, i.e.,
\[
\mathbf{d} \mathfrak{G}_c(z) \cdot v_z = 0 ~~ \forall ~~ v_z \in T_z \mathcal{C}_c^{HP}
\text{.}
\]  
Define $\mathcal{C}^Q = \{ q \in C^1([a,b], Q) ~~|~~ q(a) = i(q_1),~q(b) = i(q_2) \}$.  Let $\Psi(q)(t) = g(q(t))$ where $\Psi: \mathcal{C}^Q \to C^1([a,b],\mathbb{R}^k)$.  By the local onto theorem, there exist charts $(U, \varphi)$ of $\mathcal{C}^Q$ and a open set $V \subset \mathbb{R}^k$ such that $\Psi \circ \varphi^{-1}: U^{\prime} \times V \to V$ is a projection onto the second factor and $\varphi(q) = (x,y) \in U^{\prime} \times V$ for $q \in U$.  Set $v_q = T_{x,y} \varphi^{-1} (v_x, v_y)$ where $v_x \in U^{\prime}$ and $v_y \in V$.   Set $\tilde{\mathfrak{G}}(q) = \mathfrak{G}(q,v,p)$.  Then,
\begin{align*}
\mathbf{d} \bar{\mathfrak{G}}(z, \lambda) \cdot (v_z, v_{\lambda} ) &= 
T \tilde{\mathfrak{G}}(q) \cdot v_q + \left\langle \lambda, T \Psi(q) \cdot \frac{d}{dt} (v_q) \right\rangle \\
&= T_2 ( \tilde{\mathfrak{G}} \circ \varphi^{-1})(x,y)  \cdot  v_y  + \left\langle \lambda,    \frac{d}{dt} (v_y) \right\rangle \text{.}
\end{align*}
However, there exists $\lambda$ such that,
\[
T_2 (\tilde{\mathfrak{G}} \circ \varphi^{-1})(x,y) \cdot  v_y  + \left\langle \lambda,   \frac{d}{dt} (v_y) \right\rangle = 0 
\]
for all $v_y \in V$.   To see this expand this expression,
\[
\int_a^b \left[ \left\langle \frac{\partial L}{\partial y}(\varphi^{-1}(x,y), v), v_y \right\rangle dt +
\sum_{j=1}^m \left\langle \frac{\partial \gamma_j}{\partial y}(\varphi^{-1}(x,y)), v_y \right\rangle \circ d W_j +
\left\langle \lambda,  \frac{d}{dt} v_y \right\rangle dt  \right] = 0
\]
and note that a solution is given by:
\[
d \lambda =  \frac{\partial L}{\partial y}(\varphi^{-1}(x,y), v) dt + \sum_{j=1}^m \frac{\partial \gamma_j}{\partial y}(\varphi^{-1}(x,y)) \circ d W _j 
\]
which is an SDE for $\lambda$, with globally Lipschitz drift and diffusion vector fields.  Hence, there exists $\lambda$ such that,
\[
\mathbf{d} \bar{\mathfrak{G}}(z, \lambda) \cdot (v_z, v_{\lambda} ) = 0,~~\forall ~~ v_z \in T_z \mathcal{C}^{HP},~~v_{\lambda} \in T_{\lambda} C^0([a,b], \mathbb{R}^k)
\text{.}
\]

Mean-square symplecticity of the flows of stochastic HP equations, and hence, constrained stochastic HP equations (by the equivalence) is established in \citep{BoOw2007a}.
\end{proof}

\paragraph{Eliminating the Lagrange Multiplier}

One can eliminate the Lagrange multiplier in (\ref{eq:cstochastichp}) by taking the second Stratonovich differential of the constraint $g(q)=0$:
\[
\frac{d g}{dt}(q) = \frac{\partial g}{\partial q} \cdot v =0
\]
and,
\begin{equation}
d \frac{dg}{dt} = \frac{\partial^2 g}{\partial^2 q} ( v, v) dt + \frac{\partial g}{\partial q} \cdot dv = 0  \text{.} \label{eq:d2constraint}
\end{equation}
By differentiating the Legendre transform in (\ref{eq:cstochastichp}) we obtain,
\[
dv = \left( \frac{\partial ^2 L}{\partial ^2 v} \right)^{-1} dp - \left( \frac{\partial ^2 L}{\partial ^2 v} \right)^{-1} \left( \frac{\partial ^2 L}{\partial q \partial v} \right) v dt \text{.}
\]
Substituting $dp$ from (\ref{eq:cstochastichp}) and the above into (\ref{eq:d2constraint}) gives,
\begin{align*}
\frac{\partial^2 g}{\partial^2 q} ( v, v) dt &+ \frac{\partial g}{\partial q} \cdot  \left( \frac{\partial ^2 L}{\partial ^2 v} \right)^{-1} \left( \frac{\partial L}{\partial q} dt + \sum_{i=1}^m \frac{\partial \gamma_i} {\partial q} d W_i + \frac{\partial g}{\partial q}^* d \lambda \right) \\
&- \frac{\partial g}{\partial q} \cdot  \left( \frac{\partial ^2 L}{\partial ^2 v} \right)^{-1} \left( \frac{\partial ^2 L}{\partial q \partial v} \right) v dt  = 0 \text{.}
\end{align*}
Set $\mathbf{P}(q,v)= \frac{\partial g}{\partial q}(q)  \left( \frac{\partial ^2 L}{\partial ^2 v}(q,v) \right)^{-1} \frac{\partial g}{\partial q}(q)^*$ and solve for $d \lambda$ in the above to obtain:
\begin{align*}
d \lambda = &
- \mathbf{P}(q,v)^{-1} \frac{\partial^2 g}{\partial^2 q} ( v, v) dt - \mathbf{P}(q,v)^{-1} \frac{\partial g}{\partial q} \cdot  \left( \frac{\partial ^2 L}{\partial ^2 v} \right)^{-1} \left( \frac{\partial L}{\partial q} dt + \sum_{i=1}^m \frac{\partial \gamma_i} {\partial q} d W_i  \right) \\
& + \mathbf{P}(q,v)^{-1} \frac{\partial g}{\partial q}(q) \cdot  \left( \frac{\partial ^2 L}{\partial ^2 v} \right)^{-1} \left( \frac{\partial ^2 L}{\partial q \partial v}(q,v) \right) v dt  \text{.}
\end{align*}
Set $\mathbf{B}(q,v)=\frac{\partial g}{\partial q}(q)^* \mathbf{P}(q,v)^{-1} \frac{\partial g}{\partial q}(q) \left( \frac{\partial^2 L}{\partial^2 v} \right)^{-1}$.  Substitute $d \lambda$ into (\ref{eq:cstochastichp}) gives,
\begin{equation} \label{eq:cstochastichp2}
\begin{cases}
\begin{array}{ccl}
dq &=& v dt  \text{,} \\
d p &=&  
\left( \mathbf{Id} - \mathbf{B}(q,v) \right)
\left( \frac{\partial L}{\partial q}(q,v) dt + \sum_{j=1}^m \frac{\partial \gamma_j}{\partial q}(q) \circ d W_j \right)  \\
& &- \frac{\partial g}{\partial q}(q)^* \mathbf{P}(q,v)^{-1} \frac{\partial^2 g}{\partial^2 q}(v,v) dt 
  + \mathbf{B}(q,v) \left( \frac{\partial^2 L}{\partial q \partial v}(q,v) \right) v dt \text{,}  \\
  p &=& \frac{\partial L}{\partial v}(q,v) \text{.} 
\end{array}
\end{cases}
\end{equation}

\begin{Remark}
One can check using properties of $\mathbf{B}$ and $\mathbf{P}$, and the Stratonovich-Ito conversion formula that the Stratonovich correction term in (\ref{eq:cstochastichp2}) vanishes, and hence, (\ref{eq:cstochastichp2}) is equivalent to:
\begin{equation} \label{eq:cstochastichp3}
\begin{cases}
\begin{array}{ccl}
dq &=& v dt  \text{,} \\
d p &=&  
\left( \mathbf{Id} - \mathbf{B}(q,v) \right)
\left( \frac{\partial L}{\partial q}(q,v) dt + \sum_{j=1}^m \frac{\partial \gamma_j}{\partial q}(q) d W_j \right)  \\
& &- \frac{\partial g}{\partial q}(q)^* \mathbf{P}(q,v)^{-1} \frac{\partial^2 g}{\partial^2 q}(v,v) dt 
  + \mathbf{B}(q,v) \left( \frac{\partial^2 L}{\partial q \partial v}(q,v) \right) v dt \text{,}  \\
  p &=& \frac{\partial L}{\partial v}(q,v) \text{.} 
\end{array}
\end{cases}
\end{equation}

\end{Remark}

\section{Constrained, Stochastic VPRK Integrators}  \label{sec:stochasticvi}

\paragraph{General Case}

Let $[a,b]$ and $N$ be given, and define the fixed step size $h=(b-a)/N$ and $t_k = h k$, $k=0,...,N$.    The stochastic VPRK method will be derived  by discretizing the kinematic constraint in the stochastic HP action integral with an s-stage RK approximant (cf.~(\ref{eq:sstagerkkc})), and the integral of the Lagrangian over the kth-time step by the following quadrature:
\[
\int_{t_k}^{t_k + h} L(q,v) dt \approx h b_i L(Q_k^i, V_k^i)  \text{.}
\]  
Motivated by the methods introduced in \citep{MiReTr2003, MiTr2004}, the stochastic part of the stochastic HP action will be approximated by:
\[
\int_{t_k}^{t_k + h} \gamma_r(q) \circ d W \approx  \left( \nu_i \phi_r + \kappa_i \psi_r \right) \gamma_r(Q_k^i) \text{,}
\]
where $\phi_r: \Omega \to \mathbb{R}$ and $\psi_r: \Omega \to \mathbb{R}$ for $r=1,...,m$.  With this discretization $3/2$-order strongly convergent methods can be derived.  The constraint $g(q)=0$ is enforced for all internal stage positions $\{ Q_k^i \}_{i=1}^s$ using Lagrange multipliers as in the deterministic case.
\begin{definition}
Fix two points $q_1$ and $q_2$ on $Q$ and define the {\bfi discrete constrained, stochastic VPRK action sum} $\mathfrak{G}_d: \Omega \times \mathcal{C}_d \to \mathbb{R}$ by:
\begin{align*}
\mathfrak{G}_d =& \sum_{k=0}^{N-1} \sum_{i=1}^{s}   \left[ h b_i L(Q_k^i, V_k^i)  +
 \sum_{r=1}^m \left( \nu_i \phi_r + \kappa_i \psi_r \right) \gamma_r(Q_k^i) +h  b_i \left\langle \Lambda_k^i, g(Q_k^i) \right\rangle \right. \\
& \left. + h \left\langle p_k^i, (Q_k^i - q_k)/h -   \sum_{j=1}^s a_{ij} V_k^j \right\rangle  +  h \left\langle p_{k+1}, (q_{k+1} - q_k)/h -   \sum_{j=1}^s b_{j} V_k^j \right\rangle   \right]  \text{.}
\end{align*}
\end{definition}

In the following we assume that the VPRK method has a well-defined solution much like the assumption made in theorem~\ref{thm:cvprk}.

\begin{theorem} \label{thm:csvprk}
Given an s-stage RK method with $b$-vector and $a$-matrix, a Lagrangian system with smooth Lagrangian $L: TQ \to \mathbb{R}$ such that $\partial^2 L/ \partial v^2$ is invertible, stochastic potential $\gamma_r: Q \to \mathbb{R}$ for $r=1,...,m$, and smooth holonomic constraint function $g: Q \to \mathbb{R}^k$.     A discrete curve $c_d \in \mathcal{C}_d(q_1, q_2)$ satisfies the following stochastic VPRK method:
\begin{equation}  \label{eq:csvprk}  
\begin{cases}
\begin{array}{rll}
Q_k^i &=  & q_k + h \sum_{j=1}^s a_{ij} V_k^j,  \\
q_{k+1} &=  &q_k + h \sum_{j=1}^s b_j V_k^j,   \\
P_{k}^i &= & p_k + h \sum_{j=1}^s \left( b_j  - \frac{b_j a_{ji}}{b_i} \right) \left( \frac{\partial L}{\partial q}(Q_k^j, V_k^j) + \frac{\partial g}{\partial q}(Q_k^j)^* \cdot \Lambda_k^j \right)  \\
& & ~~~+  \sum_{j=1}^s  \sum_{r=1}^m \left(1-\frac{a_{ji}}{b_i} \right) \left( \nu_j \phi_r + \kappa_j \psi_r \right) \frac{ \partial \gamma_r}{\partial q}(Q_k^j),  \\
p_{k+1} &= & p_k + h \sum_{j=1}^s b_j \left( \frac{\partial L}{\partial q}(Q_k^j, V_k^j)  + \frac{\partial g}{\partial q}(Q_k^j)^* \cdot \Lambda_k^j \right)\\
& & ~~~+ \sum_{j=1}^s  \sum_{r=1}^m \left( \nu_j \phi_r + \kappa_j \psi_r \right) \frac{ \partial \gamma_r}{\partial q}(Q_k^j),   \\
P_k^i &= & \frac{\partial L}{\partial v} (Q_k^i, V_k^i) \text{,}  \\
g(Q^i_k) &=& 0 \text{.}
\end{array}
\end{cases}
\end{equation}
for $i=1,\cdots,s$ and $k=1,\cdots,N-1$, if and only if it is a critical point of the function $\mathfrak{G}_d:  \mathcal{C}_d \to \mathbb{R}$, that is, $\mathbf{d} \mathfrak{G}_d (c_d) = 0$.  Moreover, there exist $\{ \Lambda_k^i \}_{i=1}^s$ such that the discrete flow map defined by the above scheme, $F_h: T^*S \to T^*S$, preserves the canonical symplectic form on $T^*S$.  
\end{theorem}

\begin{proof}
We assume the existence of a numerical solution of (\ref{eq:csvprk}) and a discrete flow map $F_h:T^*S \to T^*S$.

The differential of $\mathfrak{G}_d(\omega, c_d)$  in the direction $z = (\{ \delta q_k, \delta p_k \}, \{ \delta Q_k^i, \delta V_k^i, p_k^i \}_{i=1}^s, \{ \delta \Lambda_k^i \}_{i=1}^s)$ is given by:
\begin{align*}
\mathbf{d} &\mathfrak{G}_d(\omega, c_d) \cdot z = \\
& \sum_{k=0}^{N-1} \sum_{i=1}^{s} h b_i \left[ \frac{ \partial L}{\partial q}(Q_k^i, V_k^i) \cdot \delta Q_k^i + \frac{ \partial L}{\partial v}(Q_k^i, V_k^i) \cdot \delta V_k^i  \right] + \sum_{r=1}^m \left( \nu_i \phi_r + \kappa_i \psi_r \right) \frac{ \partial \gamma_r}{\partial q}(Q_k^i) \cdot \delta Q_k^i \\
&+ h \left[ \left\langle p_k^i, (\delta Q_k^i - \delta q_k)/h -\sum_{j=1}^s a_{ij} \delta V_k^j \right\rangle  + \left\langle p_{k+1}, (\delta q_{k+1} - \delta q_k)/h -   \sum_{j=1}^s b_{j} \delta V_k^j \right\rangle \right]  \\
&+ h \left[ \left\langle \delta p_k^i, (Q_k^i - q_k)/h -\sum_{j=1}^s a_{ij} V_k^j \right\rangle  + \left\langle \delta p_{k+1}, (q_{k+1} - q_k)/h -   \sum_{j=1}^s b_{j} V_k^j \right\rangle \right]  \\
&+ h b_i \left[ \left\langle \delta \Lambda_k^i, g(Q_k^i) \right\rangle +  \left\langle \Lambda_k^i, \frac{\partial g}{\partial q}(Q_k^i) \cdot \delta Q_k^i \right\rangle \right] \text{.}
\end{align*}
Collecting terms with the same variations and summation by parts using the boundary conditions $\delta q_0 = \delta q_N = 0$ gives,
\begin{align*}
\mathbf{d} &\mathfrak{G}_d(\omega, c_d) \cdot z = \\
& \sum_{k=1}^{N-1}   \sum_{i=1}^s  \sum_{r=1}^m \left( h b_i \frac{\partial g}{\partial q}(Q_k^i)^* \Lambda_k^i  + h b_i \frac{ \partial L}{\partial q}(Q_k^i, V_k^i) +  \left( \nu_i \phi_r + \kappa_i \psi_r \right) \frac{ \partial \gamma_r}{\partial q}(Q_k^i) + p_k^i \right) \cdot \delta Q_k^i \\
&+  \left( -p_{k+1} + p_k - \sum_{i=1}^s p_k^i  \right) \cdot   \delta q_k +h \left( b_i \frac{ \partial L}{\partial v}(Q_k^i, V_k^i) - \sum_{j=1}^s a_{ji} p_k^j  - b_i p_{k+1}  \right) \cdot \delta V_k^i \\
&+ h  \left\langle \delta p_k^i, (Q_k^i - q_k)/h -\sum_{j=1}^s a_{ij} V_k^j \right\rangle  + h \left\langle \delta p_{k+1}, (q_{k+1} - q_k)/h -   \sum_{j=1}^s b_{j} V_k^j \right\rangle \\
&+ h b_i \left\langle \delta \Lambda_k^i, g(Q_k^i) \right\rangle  \text{.}
\end{align*}
Since $\mathbf{d} \mathfrak{G}_d(\omega, c_d) =0$ if and only if $\mathbf{d} \mathfrak{G}_d \cdot z = 0$ for all $z \in T_{c_d} \mathcal{C}_d$, one arrives at the desired equations with the elimination of $p_k^i$ and the introduction of the internal stage variables $P_k^i = \partial L/\partial v (Q_k^i, V_k^i) $ for $i=1,\cdots,s$.   Conversely, if $c_d$ satisfies (\ref{eq:csvprk}) then $\mathbf{d} \mathfrak{G}_d(\omega, c_d(\omega)) =0$.

The proof of symplecticity using the variational principle is very similar to the proof of theorem~\ref{thm:cvprk} and is therefore omitted.
\end{proof}

 \paragraph{Constrained, Stochastic Variational Euler}
 As an example we consider the following simple first-order strongly convergent integrator.  One can derive higher-order accurate integrators for constrained systems following the procedure for unconstrained stochastic Hamiltonian systems described in \citep{MiReTr2003, MiTr2004}.

Let $B_r^k \sim \mathcal{N}(0,h)$ be normally distributed random variables for $r=1,...,m$ and $k=0,...,N-1$.   The {\bfi constrained, stochastic variational Euler integrator} is given by:
\begin{align}
q_{k+1} =& q_k + h  \hat{v}_{k+1}  \label{eq:veulerqupdate} \text{,} \\
\hat{p}_{k+1} =& p_k + h \frac{\partial L}{\partial q}(q_k, v_k)+ \sum_{r=1}^m \frac{\partial \gamma_r}{\partial q}(q_k) B_r^k  + h \frac{\partial g}{\partial q}(q_k)^* \Lambda^1_k \text{,} \label{eq:veulerpupdate} \\
0 =& g(q_{k+1})  \text{,} \label{eq:veulerconstraint} \\
\hat{p}_{k+1} =& \frac{\partial L}{\partial v}(q_{k+1}, \hat{v}_{k+1})  \label{eq:veulerlt} \text{.} 
\end{align}  
together with the projection step (\ref{eq:projectionstep}).  Given $(q_k, p_k)$ and $h$, this integrator determines $(q_{k+1}, p_{k+1})$ by eliminating $\hat{v}_{k+1}$ using the Legendre transform, eliminating $\hat{p}_{k+1}$ using (\ref{eq:veulerpupdate}), and determining $\Lambda^1_k$ by satisfying the constraint $g(q_{k+1})=0$.  One then takes a projection step by solving (\ref{eq:projectionstep}) for $p_{k+1}$.   
The integrator defined by the composite map, $(q_k, p_k) \mapsto (q_{k+1}, \hat{p}_{k+1}) \mapsto (q_{k+1}, p_{k+1})$, has the following properties.

\begin{theorem} \label{thm:scveuler}
Constrained, stochastic variational euler is mean-squared symplectic and first-order strongly convergent integrator for (\ref{eq:cstochastichp}).
\end{theorem}

\begin{proof}
We will prove this by using the technique provided in \citep{VaCi2006}.

Consider the following first-order strongly convergent, Euler-Maruyama integrator applied to (\ref{eq:cstochastichp2}):
\begin{equation} \label{eq:firstorder}
\begin{cases}
\begin{array}{ccl}
q_{k+1} & = & q_k + h v_k  \text{,} \\
p_{k+1} & = & p_k +   \left( \mathbf{Id} - \mathbf{B}(q_k,v_k) \right)
\left( h \frac{\partial L}{\partial q}(q_k,v_k)  + \sum_{j=1}^m \frac{\partial \gamma_j}{\partial q}(q_k) B_j^k \right)  \\
& &- h \frac{\partial g}{\partial q}(q_k)^* \mathbf{P}(q_k,v_k)^{-1} \frac{\partial^2 g}{\partial^2 q}(v_k,v_k) 
  + h \mathbf{B}(q_k,v_k) \left( \frac{\partial^2 L}{\partial q \partial v}(q_k, v_k)  \cdot v_k  \right)  \text{,} \\
  p_{k+1} & = & \frac{\partial L}{\partial v}(q_{k+1}, v_{k+1})  \text{.}
\end{array}
\end{cases}
\end{equation}
The order of accuracy of constrained stochastic variational Euler will be determined by checking that the update determined by a single step of (\ref{eq:veulerqupdate})-(\ref{eq:veulerlt}) with the projection step (\ref{eq:projectionstep}) agrees with a single step of (\ref{eq:firstorder}) to $\mathcal{O}(h^2)$.

For this purpose the following expansions to $\mathcal{O}(h^2)$ are introduced: 
\begin{align*}
q_{k+1} &= q_k + h \hat{v}_{k+1},~~\hat{v}_{k+1} = v_k + h \hat{v}_{k+1}^1,~~
\hat{p}_{k+1} = p_k + h \hat{p}_{k+1}^1 \\
v_{k+1} &= v_k + h v_{k+1}^1,~~p_{k+1} = p_k + h p_{k+1}^1 
\end{align*}
Substituting these expansions into (\ref{eq:veulerqupdate}) confirms that the position update agrees to $\mathcal{O}(h^2)$.  However, checking that the momentum or velocity updates agree to $\mathcal{O}(h^2)$ is more involved.  To do this expand (\ref{eq:veulerconstraint}) in a Taylor series about $q_k$, and use $g(q_k) = 0$ and $\frac{\partial g}{\partial q}(q_k) \cdot v_k =0 $ to obtain,
\begin{equation}
g(q_{k+1}) = h^2 \frac{\partial g}{\partial q}(q_k) \cdot \hat{v}_{k+1}^1 + \frac{h^2}{2} \frac{\partial^2 g}{\partial^2 q}(q_k)(v_k, v_k) + \mathcal{O}(h^3) = 0 \label{eq:constraintexpansion}
\end{equation}
From (\ref{eq:veulerpupdate}) it is clear that,
\begin{align}
h \hat{p}_{k+1}^1 =&  h \frac{\partial L}{\partial q}(q_k,v_k)  + \sum_{j=1}^m \frac{\partial \gamma_j}{\partial q}(q_k) B_j^k + h \frac{\partial g}{\partial q}(q_k)^* \Lambda^1_k   \text{.} \label{eq:phatupdateexpansion}
\end{align}
Moreover the Taylor series expansion of (\ref{eq:veulerlt}) gives,
\begin{equation} 
h \frac{\partial^2 L}{\partial^2 v}(q_k, v_k) \cdot \hat{v}_{k+1}^1 + 
h \frac{\partial^2 L}{\partial q \partial v}(q_k, v_k) \cdot v_k = h \hat{p}_{k+1}^1 + \mathcal{O}(h^2) \label{eq:ltexpansion}
\end{equation}
Using invertibility of $\frac{\partial^2 L}{\partial^2 v}(q_k, v_k)$, one can rewrite (\ref{eq:ltexpansion}) as, 
\[
h  \hat{v}_{k+1}^1 =  - h \frac{\partial^2 L}{\partial^2 v}(q_k, v_k)^{-1}  \frac{\partial^2 L}{\partial q \partial v}(q_k, v_k) \cdot v_k + h \frac{\partial^2 L}{\partial^2 v}(q_k, v_k)^{-1} \hat{p}_{k+1}^1 + \mathcal{O}(h^2)
\]
Substitution of (\ref{eq:phatupdateexpansion}) into the above equation gives,
\begin{align*}
h  \hat{v}_{k+1}^1 &= - h \left(\frac{\partial^2 L}{\partial^2 v}\right)_k^{-1}  \left(\frac{\partial^2 L}{\partial q \partial v}\right)_k \cdot v_k \\
& +  \left(\frac{\partial^2 L}{\partial^2 v}\right)_k^{-1}\left(  h \frac{\partial L}{\partial q}(q_k,v_k)  + \sum_{j=1}^m \frac{\partial \gamma_j}{\partial q}(q_k) B_j^k + h \frac{\partial g}{\partial q}(q_k)^* \Lambda^1_k   \right) + \mathcal{O}(h^2)
\end{align*}
Substituting the above expression into (\ref{eq:constraintexpansion}) gives, 
\begin{align*}
 & -  h \left( \frac{\partial g}{\partial q} \right)_k \left(\frac{\partial^2 L}{\partial^2 v}\right)_k^{-1}  \left(\frac{\partial^2 L}{\partial q \partial v}\right)_k \cdot v_k \\
& +  \left( \frac{\partial g}{\partial q} \right)_k \left(\frac{\partial^2 L}{\partial^2 v}\right)_k^{-1}\left(  h \frac{\partial L}{\partial q}(q_k,v_k)  + \sum_{j=1}^m \frac{\partial \gamma_j}{\partial q}(q_k) B_j^k + h \frac{\partial g}{\partial q}(q_k)^* \Lambda^1_k   \right)
 \\
&+ \frac{h}{2} \frac{\partial^2 g}{\partial^2 q}(q_k)(v_k, v_k)  + \mathcal{O}(h^2) = 0
\end{align*}
From which it follows that the Lagrange multiplier satisfies, 
\begin{align*}
h \Lambda^1_k &= h \mathbf{P}_k^{-1}  \left( \frac{\partial g}{\partial q} \right)_k \left(\frac{\partial^2 L}{\partial^2 v}\right)_k^{-1}  \left(\frac{\partial^2 L}{\partial q \partial v}\right)_k \cdot v_k \\
&-  \mathbf{P}_k^{-1}  \left( \frac{\partial g}{\partial q} \right)_k \left(\frac{\partial^2 L}{\partial^2 v}\right)_k^{-1}\left(  h \frac{\partial L}{\partial q}(q_k,v_k)  + \sum_{j=1}^m \frac{\partial \gamma_j}{\partial q}(q_k) B_j^k \right)  \\
&- \frac{h}{2} \mathbf{P}_k^{-1} \frac{\partial^2 g}{\partial^2 q}(q_k)(v_k, v_k)  + \mathcal{O}(h^2) 
\end{align*}
Substitution into (\ref{eq:phatupdateexpansion}) gives, 
\begin{align}
h \hat{p}_{k+1}^1 =&  \left( \mathbf{Id} - \mathbf{B}_k \right)  \left( h \frac{\partial L}{\partial q}(q_k,v_k)  + \sum_{j=1}^m \frac{\partial \gamma_j}{\partial q}(q_k) B_j^k \right) \nonumber \\
&+ h \mathbf{B}_k   \left(\frac{\partial^2 L}{\partial q \partial v}\right)_k \cdot v_k -  \frac{h}{2} \frac{\partial g}{\partial q}(q_k)^* \mathbf{P}_k^{-1} \frac{\partial^2 g}{\partial^2 q}(q_k)(v_k, v_k)  \text{.} \label{eq:phatupdateexpansion2}
\end{align}
To determine $p_{k+1}^1$ we expand the hidden velocity constraint in (\ref{eq:projectionstep}) to obtain, 
\begin{align}
\frac{\partial g}{\partial q}(q_{k+1}) v_{k+1} & = 
h \frac{\partial g}{\partial q}(q_k) \cdot v_{k+1}^1 + h \frac{\partial^2 g}{\partial^2 q}(q_k)(v_k, v_k) + \mathcal{O}(h^2) = 0    \label{eq:hiddenvelocityexpansion}
\end{align}
Expanding the momentum update in (\ref{eq:projectionstep}) gives, 
\begin{equation} \label{eq:pupdateexpansion}
h p_{k+1}^1 = h \hat{p}_{k+1}^1 + h \frac{\partial g}{\partial q}(q_k)^* \Lambda^2_k + \mathcal{O}(h^2) \text{.}
\end{equation}
Expanding the Legendre transform in (\ref{eq:projectionstep}) gives,
\[
h v_{k+1}^1 = h \left(\frac{\partial^2 L}{\partial^2 v}\right)_k^{-1}  \hat{p}_{k+1}^1 + 
h \left(\frac{\partial^2 L}{\partial^2 v}\right)_k^{-1}  \frac{\partial g}{\partial q}(q_k)^* \Lambda^2_k - 
h \left(\frac{\partial^2 L}{\partial^2 v}\right)_k^{-1} \left( \frac{\partial^2 L}{\partial q \partial v}\right)_k \cdot v_k + \mathcal{O}(h^2)
\]
Substitution of the above into (\ref{eq:hiddenvelocityexpansion}) yields, 
\begin{align*}
 & h  \frac{\partial g}{\partial q}(q_k) \cdot  \left(\frac{\partial^2 L}{\partial^2 v}\right)_k^{-1}  \hat{p}_{k+1}^1 + 
h  \mathbf{P}_k \Lambda^2_k - 
h \frac{\partial g}{\partial q}(q_k) \cdot \left(\frac{\partial^2 L}{\partial^2 v}\right)_k^{-1} \left( \frac{\partial^2 L}{\partial q \partial v}\right)_k \cdot v_k \\
&  + h \frac{\partial^2 g}{\partial^2 q}(q_k)(v_k, v_k) + \mathcal{O}(h^2) = 0
\end{align*}
Solving the above for $\Lambda^2_k$ and substituting into (\ref{eq:pupdateexpansion}) gives,
\begin{equation*}
h p_{k+1}^1 = h (\mathbf{Id} - \mathbf{B}_k) \hat{p}_{k+1}^1 + h \mathbf{B}_k \left( \frac{\partial^2 L}{\partial q \partial v}\right)_k \cdot v_k - h \frac{\partial g}{\partial q}(q_k)^* \mathbf{P}_k^{-1}  \frac{\partial^2 g}{\partial^2 q}(q_k)(v_k, v_k)  + \mathcal{O}(h^2)
\end{equation*}
Substituting (\ref{eq:phatupdateexpansion2}) into the above and using the following identities
\[
 (\mathbf{Id} - \mathbf{B}_k) \mathbf{B}_k = 0,~~~
(\mathbf{Id} - \mathbf{B}_k)   \frac{\partial g}{\partial q}(q_k)^* \mathbf{P}_k^{-1}  \frac{\partial^2 g}{\partial^2 q}(q_k)(v_k, v_k) = 0
\] 
implies that, 
\begin{align} \label{eq:pupdateexpansion2}
p_{k+1} &= p_k +   (\mathbf{Id} - \mathbf{B}_k)  \left( h \frac{\partial L}{\partial q}(q_k,v_k)  + \sum_{j=1}^m \frac{\partial \gamma_j}{\partial q}(q_k) B_j^k \right) \nonumber \\
& + h \mathbf{B}_k \left( \frac{\partial^2 L}{\partial q \partial v}\right)_k \cdot v_k - h \frac{\partial g}{\partial q}(q_k)^* \mathbf{P}_k^{-1}  \frac{\partial^2 g}{\partial^2 q}(q_k)(v_k, v_k)  + \mathcal{O}(h^2)
\end{align}
which agrees with (\ref{eq:firstorder}) to $\mathcal{O}(h^2)$.

Mean-square symplecticity of the map $(q_k, p_k) \mapsto (q_{k+1}, \hat{p}_{k+1})$ follows from the proof of theorem~\ref{thm:csvprk}.  The projection step is also mean-square symplectic.  Therefore the composite map is mean-square symplectic with respect to the symplectic form on $T^*S$ since the space of symplectic maps forms a group.
\end{proof}

\newpage

\begin{figure}[ht]
\begin{center}
\includegraphics[scale=0.4,angle=0]{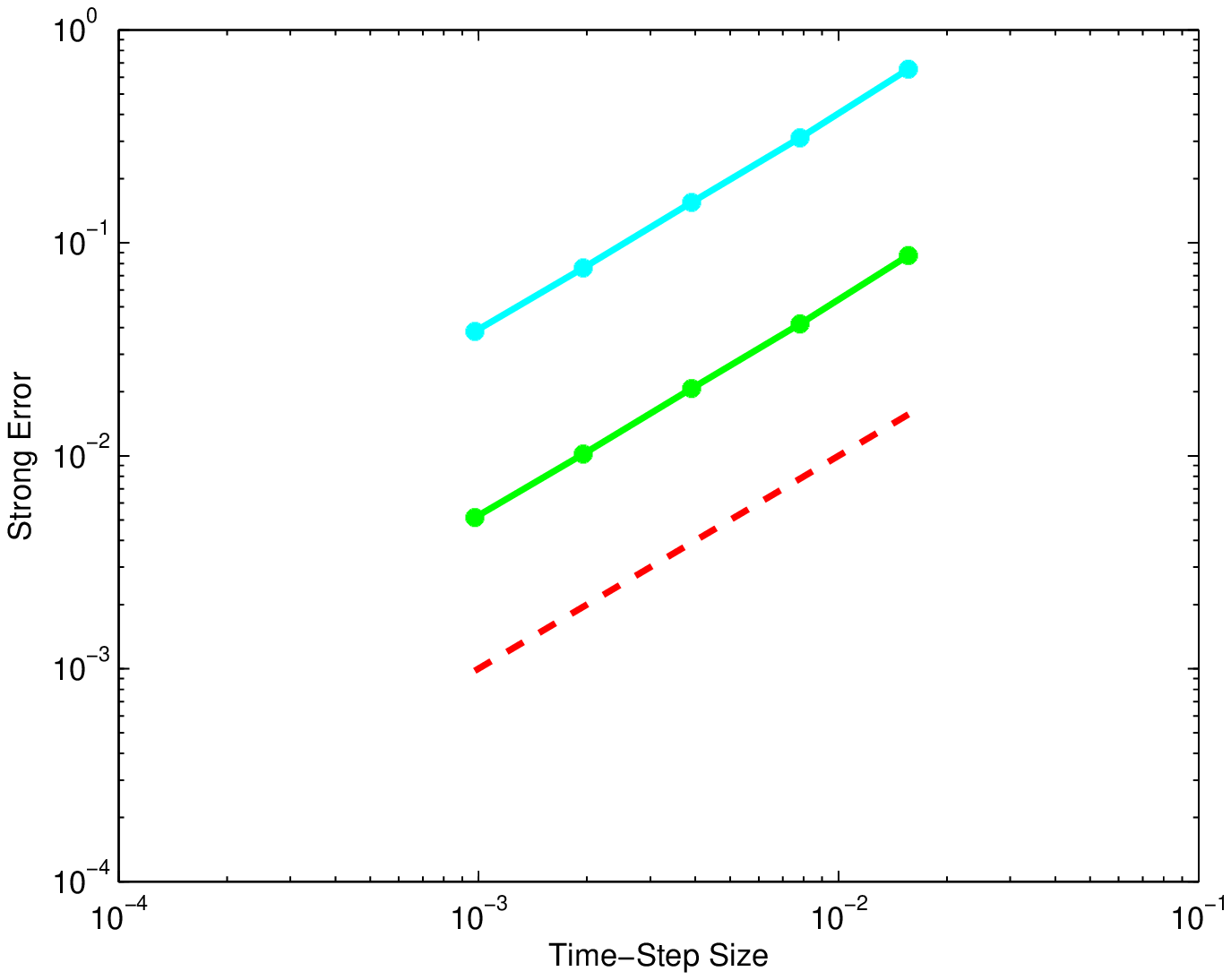}
\caption{ \small {\bf Strong error plot for stochastically perturbed spherical pendulum.}  A log-log graph of the strong error of stochastic variational Euler applied to the stochastically perturbed spherical pendulum.  Green and cyan solid lines represent the strong errors in momentum and position respectively.  For reference a dashed line of slope $1$ is included.   Observe that the order of convergence is consistent with theorem~\ref{thm:scveuler}.
}
\label{fig:strongerror}
\end{center}
\end{figure}

\begin{figure}[ht]
\begin{center}
\includegraphics[scale=0.4,angle=0]{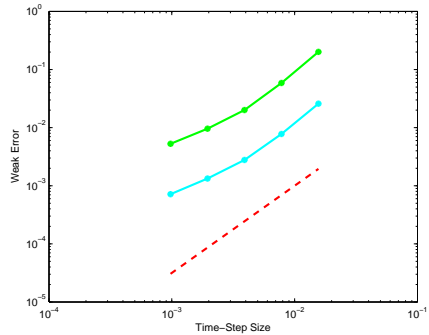}
\caption{ \small {\bf Weak error plot for stochastically perturbed spherical pendulum.} A log-log graph of the weak error of stochastic variational Euler applied to the stochastically perturbed spherical pendulum.  Green and cyan solid lines represent the strong errors in momentum and position respectively.  For reference a dashed line with slope $3/2$ is included.    This order of convergence is also consistent with theorem~\ref{thm:scveuler}.
}
\label{fig:weakerror}
\end{center}
\end{figure}

\begin{Example}[Stochastically Perturbed Spherical Pendulum]  Consider the spherical pendulum where $Q = \mathbb{R}^3$ and $S=S^2$ with holonomic constraint given by the zero level-set of $g(q) = \| q \|^2 - 1$.   Let $\{ e_1, e_2, e_3\}$ denote an orthonormal basis for $\mathbb{R}^3$ with $e_3$ corresponding to the direction of gravity.  Its Lagrangian is given by:
\[
L(q, v ) = \frac{1}{2} \| v \|^2 - q \cdot e_3 \text{.}
\] 
Consider the following stochastic potentials:  $\gamma_i(q) = \sin(q \cdot e_i)$ for $i=1,2,3$.  Applying stochastic variational Euler to the resulting stochastically perturbed spherical pendulum one obtains the strong error plot shown in Figure~\ref{fig:strongerror} confirming theorem \ref{thm:scveuler}.  Since a strong order of convergence implies at least the same weak order if not higher, it is not surprising that the weak error plot shown in Figure~\ref{fig:weakerror} seems to predict $3/2$-order weak convergence.   The dashed lines in the plots give appropriate reference slopes in each case.  
\end{Example}

\newpage

\section*{Acknowledgements}
We wish to acknowledge Jerry Marsden, Sigrid Leyendecker, and Katie Whitehead for useful suggestions on this paper.  
 
\bibliographystyle{plain}
\bibliography{nawaf}

 \end{document}